\documentclass[10pt,letterpaper,conference]{IEEEtran}
\usepackage[latin9]{inputenc}
\usepackage{amsmath}
\usepackage{amssymb}
\usepackage{graphicx}
\usepackage{esint}

\makeatletter

\usepackage{cite}

\makeatother

\begin{document}

\title{Numerical Methods for the Inverse Nonlinear Fourier Transform}

\author{\IEEEauthorblockN{Stella Civelli, Luigi Barletti}
\IEEEauthorblockA{Dipartimento di Matematica e Informatica ``U.
Dini''\\
 Università degli Studi di Firenze\\
 50134 Firenze, Italy\\
Email: stellac@garmentgroup.it, luigi.barletti@unifi.it}
\and
\IEEEauthorblockN{Marco Secondini}
\IEEEauthorblockA{TeCIP Institute\\
Scuola Superiore Sant'Anna\\
 56124 Pisa, Italy\\
Email: marco.secondini@sssup.it}}
\maketitle
\begin{abstract}
We introduce a new numerical method for the computation of the inverse
nonlinear Fourier transform and compare its computational complexity
and accuracy to those of other methods available in the literature.
For a given accuracy, the proposed method requires the lowest number
of operations.\end{abstract}

\section{Introduction}

The propagation of light in optical fibers is governed by the nonlinear
Schrödinger equation (NLSE), which accounts for the interplay of dispersion
and nonlinearity. In the lossless case, the NLSE admits an analytical
solution based on the inverse scattering transform (also known as
nonlinear Fourier transform---NFT) \cite{ZakSha72}. On this basis,
during the \textquoteright{}80s and the \textquoteright{}90s, research
on soliton transmission in optical fibers was carried out extensively.
However, despite a great effort, soliton systems never caught on commercially
and were eventually abandoned in favor of simpler and more effective
transmission techniques originally developed for  linear  channels.
Recently, an alternative approach to communication over nonlinear
channels---originally proposed in \cite{Hasegawa93} and further generalized
in \cite{Yousefi-III} and \cite{le2014nonlinear}---is gaining new
attention from the optical fiber community. Rather than using solitons
as information carriers, the idea is that of encoding information
onto the discrete and/or continuous spectrum of the Zakharov-Shabat
operator (the scattering data), which evolves trivially and linearly
along the fiber, such that the encoded information can be easily and
accurately recovered from the optical signal after propagation at
any distance, without any impact from dispersion and nonlinearity.
This approach requires the computation of the inverse NFT at the transmitter
to obtain the transmitted time-domain signal from the modulated scattering
data, and the computation of the direct NFT at the receiver to extract
the scattering data from the received time-domain signal. Both operations
are, in general, quite involved and research on efficient numerical
methods for their implementation is in progress \cite{arico2011numerical,Yousefi-II,le2014nonlinear,wahls2014fast,fermo2015scattering,wahls2015fast}.

In this work, we focus on the numerical computation of the inverse
NFT. We consider two methods available in the literature \cite{arico2011numerical,le2014nonlinear}
and propose an alternative method based on iterated convolutions evaluated
through the fast Fourier transform (FFT) algorithm. All the methods
are tested for the focusing NLSE in the soliton-free case (i.e., when
only the continuous spectrum of the Zakharov-Shabat operator is present),
which is of particular interest for information transmission based
on nonlinear inverse synthesis \cite{le2014nonlinear}. Finally, we
investigate and compare the accuracy and computational complexity
of all the methods.

\section{Numerical methods}

We consider the focusing NLSE in the normalized form 
\begin{equation}
ju_{z}+\frac{1}{2}u_{tt}+u|u|^{2}=0\label{eq:nlse}
\end{equation}
We assume to know the scattering data at a given distance $z$ and
we focus on the inverse scattering problem (or inverse NFT) to obtain
the corresponding NLSE solution $u(t,z)$.%
\footnote{Hereafter, the dependence on $z$ is inessential and will be omitted. %
} This problem leads to the Gel'fand-Levitan-Marchenko equation (GLME)
\cite{ablowitz1981solitons} 
\begin{equation}
\! K(t,y)+F(t+y)+\int\limits _{-\infty}^{t}\!\!\int\limits _{-\infty}^{t}\! K(t,r)F^{*}(r+s)F(s+y)dsdr=0\label{eq:glme}
\end{equation}
where the integral kernel $F(y)$ is a function of the scattering
data. Specifically, considering the soliton-free case, $F(y)=\frac{1}{2\pi}\int_{-\infty}^{+\infty}r(\lambda)e^{-j\lambda y}d\lambda$
where $r(\lambda)$ is the reflection coefficient. The solution of
the NLSE is finally obtained as $u(t)=2\lim_{y\rightarrow t^{-}}K(t,y)$.

Let $F'(y)=F(y)$ for $y\ge0$ and $F'(y)=0$ otherwise. The GLME
can be rewritten in the form of the Marchenko integral equations 
\begin{equation}
\begin{array}{l}
\!\!\! B_{1}(t,\alpha)=-\int\limits _{0}^{+\infty}\!\! F'^{*}(2t\!-\!\alpha\!-\!\beta)B_{2}(t,\beta)d\beta\\
\!\!\! B_{2}(t,\alpha)=\!\!\int\limits _{0}^{+\infty}\!\! F'(2t\!-\!\alpha\!-\!\beta)B_{1}(t,\beta)d\beta-F'(2t\!-\!\alpha)
\end{array}\label{eq:marchenko}
\end{equation}
and the NLSE solution obtained as $u(t)=2B_{2}(t,0^{+})$. In the
following we assume to know $F$ on a uniform grid in the interval
$[0,2T]$ with step $\delta_{\alpha}=\frac{2T}{N_{F}}$ and we analyze
three different numerical methods to solve the system \eqref{eq:marchenko},
in order to find the solution of \eqref{eq:nlse} for $t\in[0,T]$.
Specifically, we consider the uniform grid on $[0,T]$ given by $t_{m}=(m-1)\delta_{t}$
for $m=1,..,N_{u}+1$ with discretization step $\delta_{t}=\frac{T}{N_{u}}$.
The complexity of the methods is measured by considering the required
number of complex products, assuming that the FFT of a vector of $N$
elements requires $\frac{N}{2}\log_{2}(N)$ complex products.

An efficient numerical way to compute the inverse NFT, assuming $N_{F}=N_{u}$,
and therefore $\delta_{t}=\frac{\delta_{\alpha}}{2}$, is presented
in \cite{le2014nonlinear}. The authors, using the Nyström method
with a rectangular quadrature scheme to approximate the integrals
in (\ref{eq:marchenko}), find the solution $u(t_{m})$ by solving
a linear system characterized by a Toeplitz matrix of size $N_{F}+m$.
The method of Trench for the inversion of non-Hermitian Toeplitz matrices
\cite{zohar1969toeplitz} is used to recursively compute the solution
for every value $t_{m}$. For this reason, let us name this method
Nyström-Trench (NT). The resulting algorithm computes the solution
$u(t)$ 
in the desired interval with complexity%
\footnote{The computational complexity of the methods presented in \cite{le2014nonlinear}
and \cite{arico2011numerical} is not indicated by the authors. The
values reported here refer to the most efficient implementations that
we were able to devise.\label{fn:The-computational-complexity}%
} 
\begin{equation}
C\simeq11N_{F}^{2}\label{eq:Complexity_NT}
\end{equation}
However, due to the recursive nature of the algorithm, the complexity
remains almost unchanged if the solution is computed in a lower number
of points $N_{u}$.

Another technique to solve system \eqref{eq:marchenko} with $\delta_{t}=c\frac{\delta_{\alpha}}{2}$
for $c\in\mathbb{N}$, (i.e., in $N_{u}=N_{F}/c$ points in the interval
$[0,T]$) is presented in \cite{arico2011numerical}, where the Nyström
method with composite Simpson's quadrature rule is considered. For
a fixed time $t_{m}$, the resulting linear system is solved by means
of the conjugate gradient method. Therefore, we refer to this scheme
as the Nyström conjugate gradient (NCG) method. Since we deal only
with diagonal and upper-left triangular Hankel matrices, the complexity
of the method can be reduced by using the FFT algorithm to compute
the matrix-vector products required by the conjugate gradient method.
As a result, the number of operations required to compute the solution
in the desired interval is\textsuperscript{\ref{fn:The-computational-complexity}}
\begin{equation}
\begin{split}C\simeq & (1+k_{\text{max}})\sum_{m=1}^{N_{u}+1}6[c(m-1)+1]\log_{2}\{2[c(m-1)+1]\}+\\
 & +N_{u}N_{F}(4+6k_{\text{max}})
\end{split}
\label{eq:Complexity_NCG}
\end{equation}
where $k_{\mathrm{max}}$ is the number of iterations of the conjugate
gradient method, which is assumed to be the same for every $t$ for
the sake of simplicity. In this case, the solution is independently
evaluated at each time $t_{m}$ in the given interval, such that $C$
depends significantly on the  number of points $N_{u}$.

Now, let us define a new iterative method to solve the system \eqref{eq:marchenko}
for $t\in[0,T]$ with discretization step $\delta_{t}=c\frac{\delta_{\alpha}}{2}$.
The main idea is that, by defining the auxiliary functions $B_{i}'(t,\alpha)=B_{i}(t,\alpha)$
for $\alpha\ge0$ and $B_{i}'(t,\alpha)=0$ otherwise, with $i=1,2$,
the integrals in \eqref{eq:marchenko} can be seen as convolutions
between these functions and the kernel $F$ and efficiently computed
by means of the FFT algorithm. Therefore, considering a proper starting
value for $B_{1}$ and iteratively updating the values of $B_{2}$
and $B_{1}$ by alternately computing the integrals in \eqref{eq:marchenko},
we define an iterative method that, under certain conditions on $F$
and $t$, converges to the solution of the system. As a starting point
for $B_{1}$ at $t_{m}$, we pick the solution found at the previous
point $t_{m-1}$, as we are dealing with continuous functions. The
resulting algorithm is 
\begin{align}
\begin{array}{l}
\!\text{For }k=0\\
\;\; B_{1}^{(0)}(\alpha,t_{m})=\left\{ \begin{array}{ll}
0 & \text{if }m=1\\
B_{1}^{(k_{\text{max}})}(\alpha,t_{m-1}) & \text{if }m>1
\end{array}\right.\\
\!\text{For }k=1,..,k_{\text{max}}\\
\;\;\left\{ \!\!\begin{array}{l}
B_{2}^{(k)}(\alpha,t_{m})=-F(2t_{m}\!-\alpha)+(F\star{B'}_{1}^{(k-1)})(2t_{m}\!-\alpha)\\
B_{1}^{(k)}(\alpha,t_{m})=-(F^{*}\star{B'}_{2}^{(k)})(2t_{m}-\alpha)
\end{array}\right.
\end{array}\label{eq:ICmethod}
\end{align}
where $\star$ denotes convolution, $B_{i}^{(k)}$ is the $k$th estimate
of the function $B_{i}$, and ${B'}{}_{1}^{(k)}$ the corresponding
auxiliary function. We refer to this method as the iterative convolution
(IC) method. The number of operations required to find the solution
in the $N_{u}$ points inside the interval $[0,T]$ is 
\begin{equation}
C\simeq k_{\text{max}}\!\!\sum_{m=1}^{N_{u}+1}\!3[2c(m-1)+1]\log_{2}[2c(m-1)+1]+2k_{\text{max}}N_{u}N_{F}\label{eq:Complexity_IC}
\end{equation}
Also in this case, $C$ depends significantly on the number of points
$N_{u}$ and can be reduced by evaluating the solution with a larger
step size.

\section{Numerical results}

With the purpose of comparing the methods presented above, we consider
the single-pulse signal considered in \cite{le2014nonlinear} 
\begin{equation}
u(t)=-\frac{4\alpha\nu\sigma(\sigma-1)}{(\sigma-1)^{2}e^{-2\sigma\alpha t}+\nu^{2}e^{2\sigma\alpha t}}\label{eq:segn}
\end{equation}
where $\sigma=\sqrt{\nu^{2}+1}$. In this case, the scattering data
are known and the GLME kernel is \cite{le2014nonlinear} 
\begin{equation}
F(y)=\alpha\nu e^{-\alpha y}
\end{equation}
for $\alpha>0$ and $-1\leq\nu\leq1$. We consider the interval $[0,3]$
and show the results only for $\alpha=\nu=1$, as similar results
are obtained for different values of the parameters.

\begin{figure}
\centering{}\includegraphics[width=0.77\columnwidth]{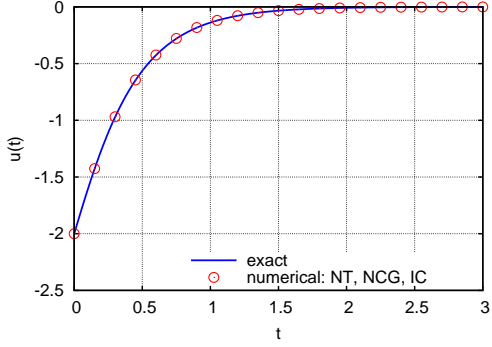}\vspace*{-0.8ex}
\caption{\label{fig:pulse_u}Comparison between the exact solution \eqref{eq:segn}
and the solution obtained with the NT, NCG, and IC numerical methods.}
\end{figure}
\begin{figure}
\begin{centering}
\vspace*{-4.5ex}
\hspace*{5ex}\includegraphics[bb=60bp 55bp 385bp 300bp,width=0.5\columnwidth]{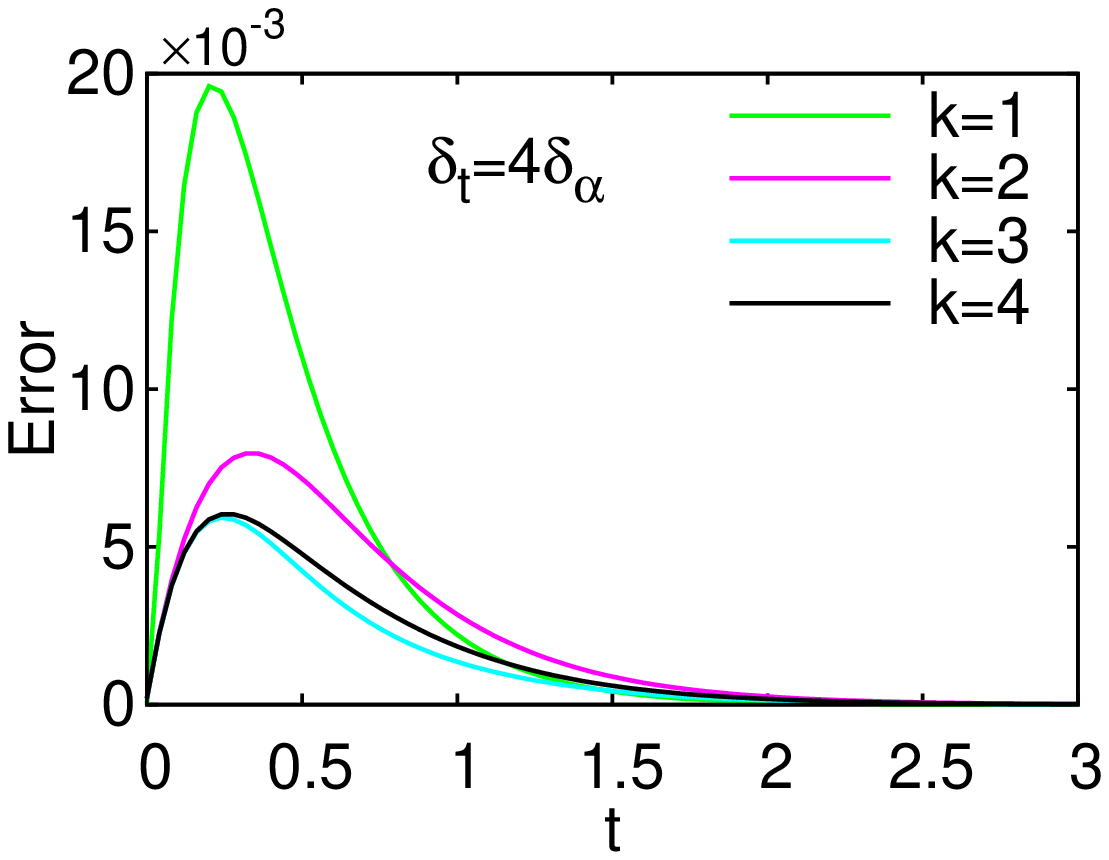}\includegraphics[bb=60bp 55bp 385bp 300bp,width=0.5\columnwidth]{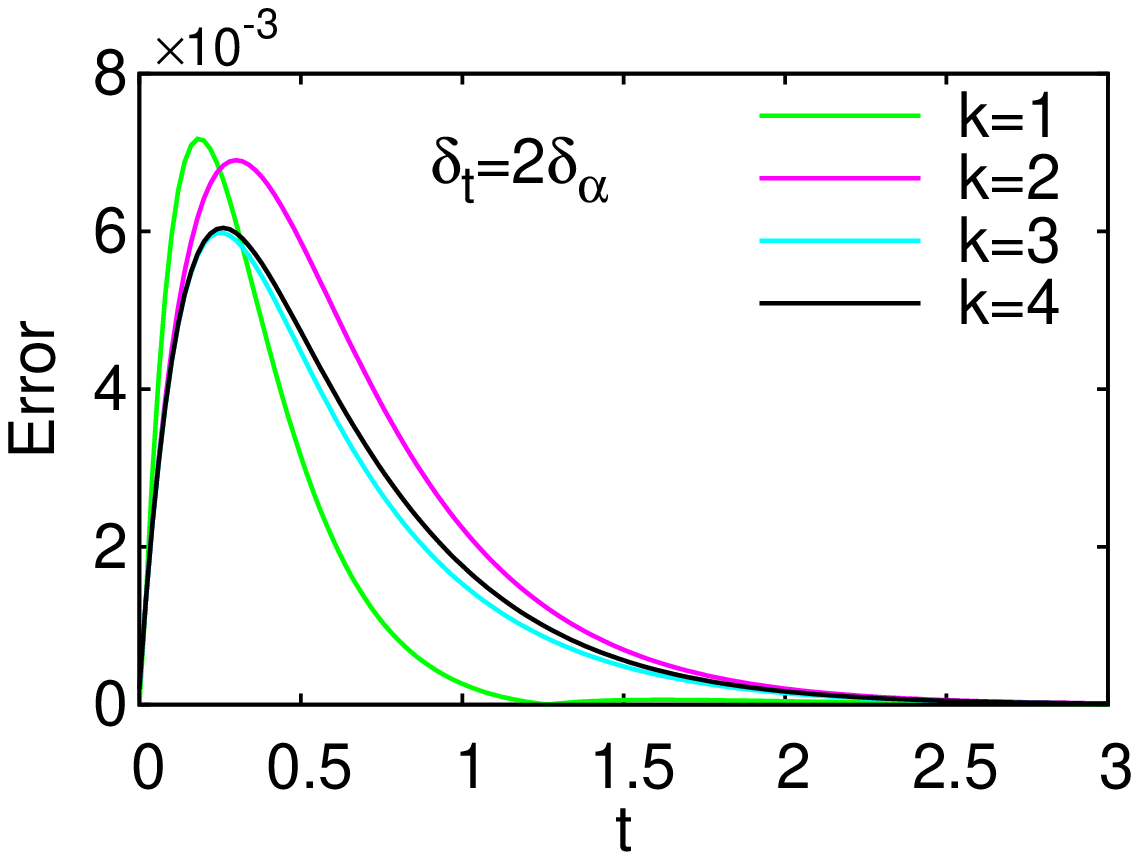}
\par\end{centering}

\centering{}\hspace*{5ex}\includegraphics[bb=60bp 55bp 385bp 300bp,width=0.5\columnwidth]{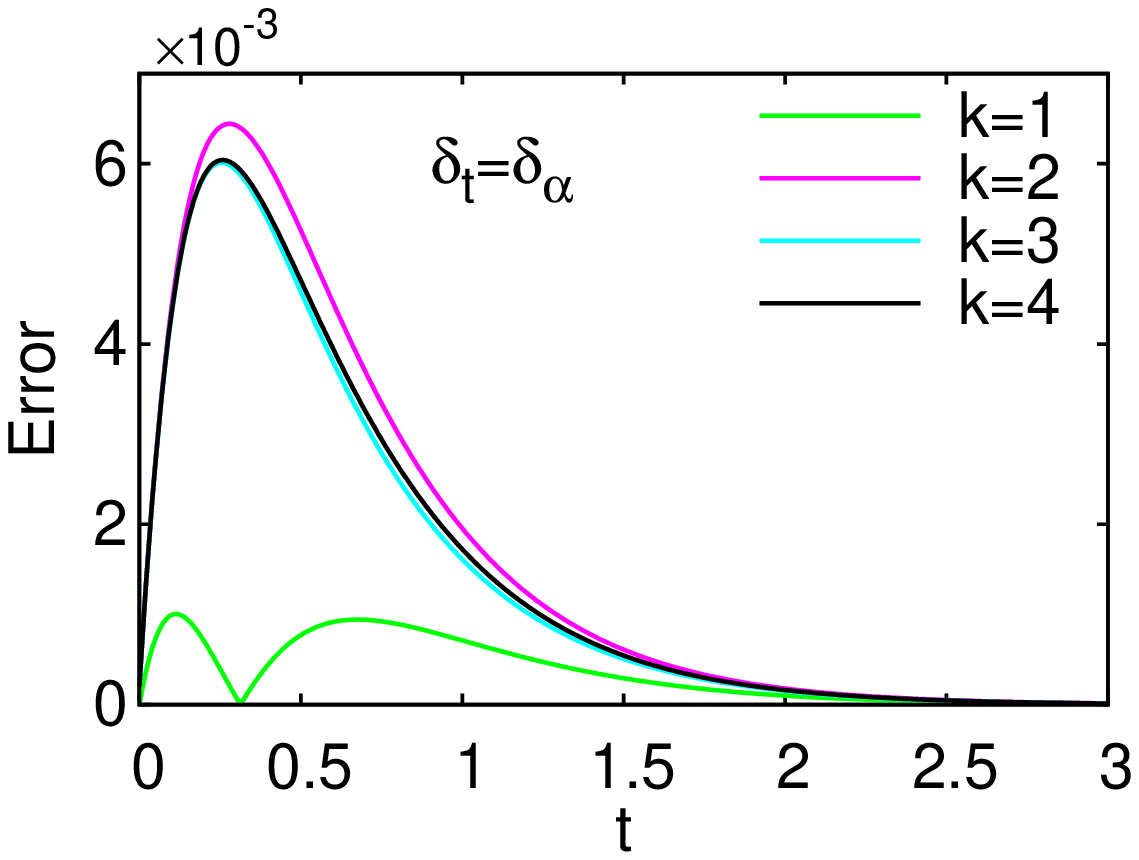}\includegraphics[bb=60bp 55bp 385bp 300bp,width=0.5\columnwidth]{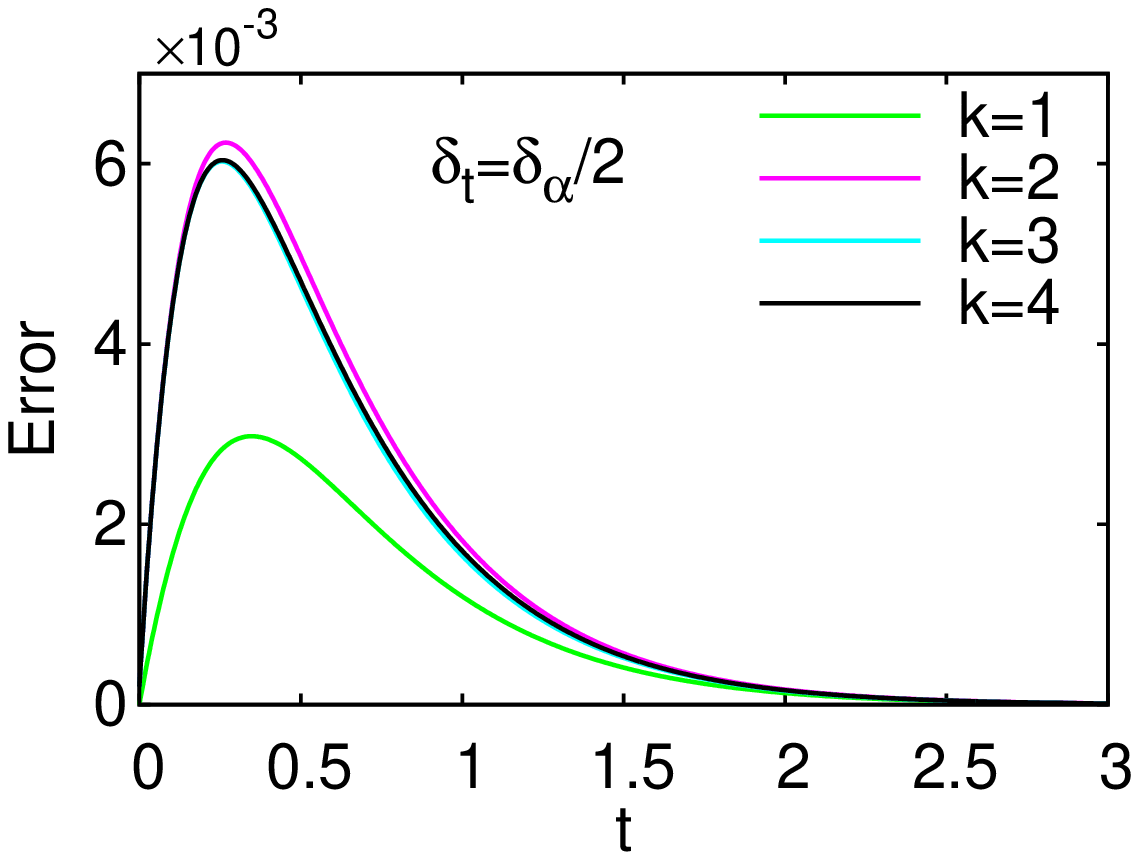}\vspace*{4ex}
\caption{\label{fig:error-vs-deltau}IC method: modulus of the error for different
iterations $k$, different resolutions $\delta_{t}$, and a fixed
resolution $\delta_{\alpha}=0.01$.}
\end{figure}

In Fig.~\ref{fig:pulse_u} we compare the analytical solution \eqref{eq:segn}
with the numerical results obtained with $\delta_{\alpha}=0.01$ and
$\delta_{t}=0.005$. At this scale, all the numerical results are
practically superimposed (and are, therefore, represented by a single
curve to avoid cluttering the figure) and equal the exact solution.
Since the new IC method is an iterative method, whose starting point
at $t_{m}$ depends on the previous solution found at $t_{m-1}$,
it is interesting to investigate its convergence properties when the
solution is evaluated with a different resolution $\delta_{t}$ (keeping
unchanged the resolution $\delta_{\alpha}$ for the evaluation of
the integrals). Fig.~\ref{fig:error-vs-deltau} shows the modulus
of the error, i.e., the difference between (\ref{eq:segn}) and the
numerical approximation, as a function of $t$ for $\delta_{\alpha}=0.01$,
different iteration numbers $k$, and different resolutions $\delta_{t}$.
In all the cases, three iterations are enough to reach the convergence.
However, in the case $\delta_{t}=\delta_{\alpha}$, the error after
the first iteration is significantly lower than the error obtained
after the convergence is reached. This quite surprising result has
been observed also for different resolutions $\delta_{\alpha}$ and
for different pulse shapes (not shown here due to a lack of space).
However, we have not been able to find a good explanation for this
behavior, which is still under investigation. Given its particular
properties, in the rest of the paper we will refer to the IC method
with $k_{\mathrm{max}}=1$ and $\delta_{t}=\delta_{\alpha}$ as the
IC1 method.

\begin{figure}
\centering{}\includegraphics[width=0.77\columnwidth]{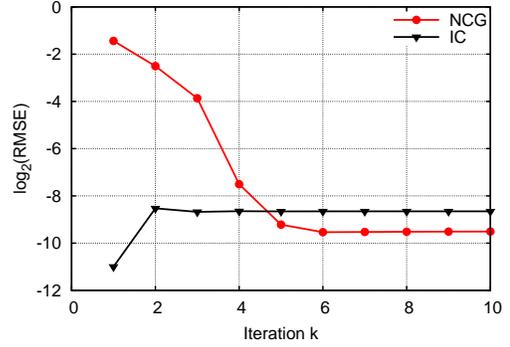}\vspace*{-0.8ex}
 \caption{\label{fig:RMSE-vs-iteration}RMSE versus number of iterations for
the NCG and the IC method with $\delta_{\alpha}=0.01$ (and $\delta_{t}=\delta_{\alpha}$
for IC).}
\end{figure}

Fig.~\ref{fig:RMSE-vs-iteration} compares the convergence properties
of the IC method (for $\delta_{t}=\delta_{\alpha}=0.01$) to those
of the NCG method (for $\delta_{\alpha}=0.01$ and any $\delta_{t}$,
as they are independent of $\delta_{t}$), showing how the root-mean-square
error (RMSE) evaluated over the whole interval $[0,3]$ evolves with
the number of iterations. 
Thanks to the composite Simpson's quadrature rule employed by the
NCG method (which is more accurate than the rectangular quadrature
rule implicitly employed by the IC method when using the FFT to compute
the convolutions), the final RMSE reached by the NCG method is lower
than the RMSE reached by the IC method. However, the IC method converges
faster (after 3 iterations) than the NCG method (after 6 iterations).
The same fast convergence of the IC method is observed also for higher
values of $\delta_{t}$ (as shown in Fig.~\ref{fig:error-vs-deltau}).
However, in the special case $\delta_{t}=\delta_{\alpha}$, as observed
before, the IC method achieves an even lower RMSE at the first iteration
(IC1 method). Such a fast convergence is obtained thanks to the good
starting value for $B_{1}$ indicated in (\ref{eq:ICmethod}). A
convergence of the iterative method is observed in all the considered
cases. This is, however, not guaranteed in general, as the convergence
depends on specific conditions on $F$ and $t$ that are currently
being investigated.

\begin{figure}
\centering{}\includegraphics[width=0.77\columnwidth]{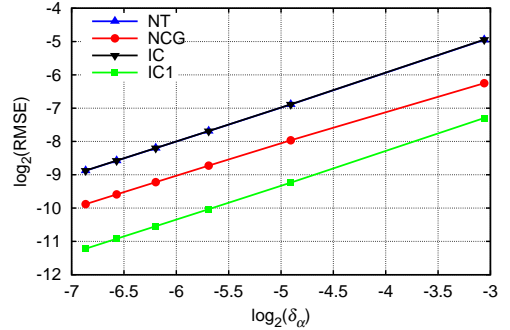}\vspace*{-0.8ex}
 \caption{\label{fig:RMSE-vs-deltaF}RMSE as a function of the resolution $\delta_{\alpha}$
for all the numerical methods ($\delta_{t}=\delta_{\alpha}/2$ for
NT and $\delta_{t}=\delta_{\alpha}$ for IC1).}
\end{figure}

For all the methods, both the accuracy and the complexity depend on
the resolution $\delta_{\alpha}$ employed for the quadrature rules.
Fig.~\ref{fig:RMSE-vs-deltaF} shows the RMSE as a function of $\delta_{\alpha}$
for all the considered methods. 
In the NT and IC1 methods, the resolution $\delta_{t}$ is set to
$\delta_{\alpha}/2$ and $\delta_{\alpha}$ by construction, respectively.
On the other hand, in the NCG and IC methods, $\delta_{t}$ can be
selected as an integer multiple of $\delta_{\alpha}/2$ without affecting
the accuracy. The RMSE values obtained with the NT and IC methods
are practically the same, as they are both based on a rectangular
quadrature rule for the computation of the integrals in \eqref{eq:marchenko}.
The NCG method, thanks to the more accurate composite Simpson's quadrature
rule, achieves a lower RMSE for the same resolution $\delta_{\alpha}$.
Finally, the IC1 method achieves an even lower RMSE.

\begin{figure}
\centering{}\includegraphics[width=0.77\columnwidth]{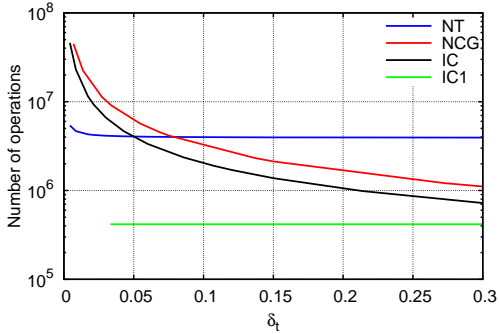}\vspace*{-0.8ex}
 \caption{\label{fig:complexity-vs-deltau}Number of operations required by
the four methods for a fixed RMSE of $2\times10^{-3}$ and a variable
desired resolution $\delta_{t}$.}
\end{figure}

Fig.~\ref{fig:RMSE-vs-deltaF} can be used to determine, for each
method, what is the resolution $\delta_{\alpha}$ (and, correspondingly,
the number of points $N_{F}=2T/\delta_{\alpha}$) required for the
computation of the integrals in order to find the solution $u(t)$
with a given accuracy. This, in turn, will affect the complexity of
the method according to \eqref{eq:Complexity_NT}, \eqref{eq:Complexity_NCG},
or \eqref{eq:Complexity_IC}. Another relevant parameter that affects
the computational complexity of the NCG and IC methods (but not that
of the NT and IC1 methods) is the resolution $\delta_{t}=c\delta_{\alpha}/2$
with which the solution $u(t)$ is needed (or, equivalently, the number
of points $N_{u}=N_{F}/c$ in which the solution must be computed).
Fig.~\ref{fig:complexity-vs-deltau} compares the computational complexity
required by the various methods, computed according to (\ref{eq:Complexity_NT}),
(\ref{eq:Complexity_NCG}), and (\ref{eq:Complexity_IC}), to obtain
the solution $u(t)$ with a desired RMSE of $2\times10^{-3}$ (about
$2^{-9}$ in Fig.~\ref{fig:RMSE-vs-deltaF}), as a function of the
desired resolution $\delta_{t}$. 
 According to Fig.~\ref{fig:RMSE-vs-deltaF}, the required resolution
is $\delta_{\alpha}\simeq0.0086$ for the NT and IC methods, $\delta_{\alpha}\simeq0.014$
for the NCG method, and $\delta_{\alpha}\simeq0.033$ for the IC1
method. The number of iterations required to reach the convergence
is $k_{\mathrm{max}}=6$ for the NCG method and $k_{\mathrm{max}}=3$
for the IC method. By construction, the minimum resolution $\delta_{t}$
is $\delta_{\alpha}/2$ for the NT, NCG, and IC methods, and $\delta_{\alpha}$
for the IC1 method. However, if desired, the resolution $\delta_{t}$
can be increased by computing the solution $u(t)$ in a lower number
of points $N_{u}$ in the given interval $[0,T]$. This will significantly
reduce the computational complexity of the NCG and IC methods, but
will leave the complexity of the NT and IC1 methods almost unaffected.%
\footnote{Even if the solution is needed in a lower number of points $N_{u}<N_{F}$,
the NT and IC1 algorithms must be anyway executed for $N_{F}$ points.
In the NT case, however, some operations can be saved, slightly reducing
the computational complexity with respect to what is reported in \eqref{eq:Complexity_NT}.%
} For any resolution $\delta_{t}$, the NCG and IC methods have a similar
complexity. In particular, the latter is slightly more convenient,
though the former could be probably improved by using a more accurate
quadrature rule or employing preconditioning. Both methods are more
complex than the NT method when employed at full resolution, but become
less complex when the solution is required in a lower number of points
(the IC method for $\delta_{t}>0.05$ and the NCG method for $\delta_{t}>0.075$).
Finally, the IC1 method is significantly less complex than the NT
method (by about one order of magnitude) and less complex than the
NCG and IC methods even at high resolutions $\delta_{t}$. However,
the behavior of the IC1 method is still under investigation and its
superior performance, though observed also in other cases (not reported
in this paper for a lack of space), has not been fully explained nor
demonstrated for a generic signal.

\section{Conclusions}

A new numerical method for the computation of the inverse nonlinear
Fourier transform has been proposed. The solution is obtained after
iterated convolutions with the GLME kernel, which are efficiently
computed trough the FFT algorithm. The accuracy and computational
complexity of the proposed method have been investigated and compared
to those of two other methods available in the literature, both exploiting
the Nyström method in combination with either the conjugate gradient
algorithm \cite{arico2011numerical} or the Trench algorithm \cite{le2014nonlinear}.
The obtained results show that, for a desired accuracy, the proposed
method requires the lowest number of operations. The relation between
the discretization step with which the integrals of the GLME are approximated
and the time resolution with which the solution is computed (which
need not be the same) is also discussed for all the methods, and the
dependence of accuracy and complexity on those two parameters is investigated.
When the two parameters are equal, the proposed method shows a very
good but peculiar behavior which is still under investigation. Further
investigation is also required to compare the algorithms in different
scenarios, such as for the defocusing NLSE or in the presence of discrete
eigenvalues in the nonlinear spectrum of the signal.

\section*{Acknowledgment}

This work was supported in part by the Italian MIUR under the FIRB
project COTONE

\bibliographystyle{IEEEtran}
\bibliography{ref}

\end{document}